\documentclass[11pt]{article}
\usepackage{amsmath}

\usepackage{amsfonts}
\usepackage{amssymb}

\newlength{\minpag}
\input{epsf}
\input epsf.tex
\begin{document}

\title{{\bf \Large Harmonic Univalent Mappings and Linearly Connected Domains }}

\author{Martin Chuaqui  \and Rodrigo Hern\'andez \thanks{The
authors were partially supported by Fondecyt Grant \# 1030589.
\endgraf \hspace{4pt} {\sl Key words: Harmonic mapping, univalent, linearly connected,
second complex dilatation.}
\endgraf \hspace{4pt} {\sl 2000 AMS Subject Classification}
$\#$: 30C99, 31A05.
\endgraf
\hspace{4pt} {\sl Email:} mchuaqui@mat.puc.cl;
rodrigo.hernandez@uai.cl.} }

\date{}
\maketitle

\begin{abstract}
We investigate the relationship between the univalence of $f$ and
of $h$ in the decomposition $f=h+\overline{g}$ of a
sense-preserving harmonic mapping defined in the unit disk
$\mathbb{D}\subset\mathbb{C}$. Among other results, we determine
the holomorphic univalent maps $h$ for which there exists $c>0$
such that every harmonic mapping of the form $f=h+\overline{g}$
with $|g'|< c|h'|$ is univalent. The notion of a linearly
connected domain appears in our study in a relevant way.
\end{abstract}

A planar harmonic mapping is a complex-valued harmonic function
$f(z)$, $z=x+iy$, defined on some domain
$\Omega\subset\mathbb{C}$. When $\Omega$ is simply connected, the
mapping has a canonical decomposition $f=h+\overline{g}$, where
$h$ and $g$ are analytic in $\Omega$.  Since the Jacobian of $f$
is given by $|h'|^2 - |g'|^2$, it is locally univalent and
orientation-preserving if and only $|g'|< |h'|$, or equivalently
if $h'(z)\neq0$ and the dilatation $\omega=g'/h'$ has the property
$|\omega(z)|<1$ in $\Omega$.

Fundamental questions regarding univalent harmonic mappings are
still to be resolved, including important coefficient estimates,
and the exact nature of the analogue of the Riemann mapping
theorem. There are beautiful results such as the shear
construction of Clunie and Sheil-Small [C-SS], and the theorem of
Rad\'o-Kneser-Choquet for convex harmonic mappings [K], [Ch]. The
literature nevertheless appears to contain few results about such
basic issues as the relation between the univalence of $f$ and of
$h$. In this paper we determine conditions under which the
univalence of one of them implies that of the other. We also find
conditions under which the harmonic mappings
$F=h+e^{i\theta}\overline{g}$ remain univalent for all
$\theta\in[0,2\pi]$.

 A domain $\Omega\subset\mathbb{C}$ is {\it linearly
connected} if there exists a constant $M<\infty$ such that any two
points $w_1,w_2\in\Omega$ are joined by a path
$\gamma\subset\Omega$ of length $l(\gamma)\leq M|w_1-w_2|$, or
equivalently (see [P]), ${\rm diam}(\gamma)\leq M|w_1-w_2|$. Such
a domain is necessarily a Jordan domain, and for piecewise
smoothly bounded domains, linear connectivity is equivalent to the
boundary's having no inward-pointing cusps. Our first result can
be considered a harmonic version of [A-B].

\medskip
\noindent{\bf Theorem 1:} {\sl Let
$h:\mathbb{D}\rightarrow\mathbb{C}$ be a holomorphic univalent
map. Then there exits $c>0$ such that every harmonic mapping
$f=h+\overline{g}$ with dilatation $|\omega|<c$ is univalent if
and only if $h(\mathbb{D})$ is a linearly connected domain.}

\medskip

The proof will show that $c$ may be taken equal to 1 when $h$ is
convex, and we will show that $c=1$ only in this case.

\medskip
\noindent{\bf Proof:} Suppose first that $\Omega=h(\mathbb{D})$ is
linearly connected, and let $f=h+\overline{g}$ be a harmonic
mapping with $|\omega|<1/M$. We claim that $f$ is univalent in
$\mathbb{D}$, or equivalently, that $w+\overline{\varphi(w)}$ is
univalent in $\Omega$, where $\varphi=g\circ h^{-1}$ satisfies
$|\varphi'|<1/M$. If not, then for $w_1\neq w_2$ we will have
$$ \varphi(w_2)-\varphi(w_1) = \overline{w_1-w_2} \, . \eqno (1)$$
Let $\gamma\subset \Omega$ be a path joining $w_1$ and $w_2$ with
$l(\gamma)\leq M|w_1-w_2|$. Then
$$ |\varphi(w_2)-\varphi(w_1)| \leq \int_{\gamma}|\varphi'(w)||dw|
< \frac{1}{M}\int_{\gamma}|dw| \leq |w_1-w_2|\, , $$ which
contradicts (1). When $\Omega$ is convex then we may take $M=1$,
and thus $f$ will be univalent as long as $|\omega|< 1$.

\smallskip
Suppose now that $\Omega$ is not linearly connected. We seek
nonunivalent harmonic mappings of the form $f=h+\overline{g}$ with
arbitrarily small dilatation. According to the main result in
[A-B], for every $c>0$ there exists a nonunivalent holomorphic
function $\phi:\Omega\rightarrow\mathbb{C}$ that satisfies
$|\phi'(w)-1|<c\,$ on $\Omega$. We write $\phi(w)=w+\varphi(w)$,
so that $|\varphi'(w)|<c$. Let $w_1,w_2$ be distinct points for
which $\phi(w_1)=\phi(w_2)$, that is, points for which
$w_1-w_2=\varphi(w_2)-\varphi(w_1)$. Let $\theta$ be such that
$e^{i\theta}(w_1-w_2)\in\mathbb{R}$. We claim that the harmonic
mapping $F(w)=w+e^{-2i\theta}\overline{\varphi(w)}$ fails to be
univalent in $\Omega$. Indeed,
$$w_1-w_2=e^{-i\theta}e^{i\theta}(w_1-w_2)
=e^{-i\theta}\overline{e^{i\theta}(w_1-w_2)}=e^{-2i\theta}
\overline{w_1-w_2}=e^{-2i\theta}
\overline{\varphi(w_2)-\varphi(w_1)}\, ,$$ that is,
$F(w_1)=F(w_2)$. The nonunivalent harmonic mapping $f(z)=F(h(z))$
has dilatation bounded by $c$. This finishes the proof.

\bigskip

We point out that in the case that $\Omega=h(\mathbb{D})$ is
linearly connected with constant $M$ and $|\omega|\leq c$ for some
$c<1/M$, then $R=f(\mathbb{D})$ is also linearly connected. In
effect, let us regard $R$ as the image of $\Omega$ under the
harmonic mapping $\phi(w)=w+\overline{\varphi(w)}$, where
$\varphi=g\circ h^{-1}$ satisfies $|\varphi'|\leq c$. Let
$\zeta_1=\phi(w_1),\, \zeta_2=\phi(w_2)\in R$ be distinct points.
By assumption, there exists $\gamma\subset\Omega$ joining $w_1$
and $w_2$ such that $l(\gamma)\leq M|w_1-w_2|$. Let
$\Gamma=\phi(\gamma)$. Because the complex derivatives of $\phi$
satisfy $|\phi_w|+|\phi_{\overline{w}}|\leq 1+c$, it follows that
$|l(\Gamma)|\leq (1+c)l(\gamma)\leq (1+c)M|w_1-w_2|$. On the other
hand,
$$|\zeta_1-\zeta_2|=|w_1-w_2+\overline{\varphi(w_1)-\varphi(w_2)}|
\geq \,|w_1-w_2|-|\varphi(w_1)-\varphi(w_2)|
$$
$$\hspace{3.25cm} \geq \, |w_1-w_2|-\int_{\gamma}|\varphi'(w)||dw| \geq
|w_1-w_2|-cl(\gamma) \geq (1-cM)|w_1-w_2| \, , $$ and therefore
$$l(\Gamma)\leq \frac{(1+c)M}{1-cM}\, |\zeta_1-\zeta_2| \, ,$$
that is, $R$ is linearly connected with constant $(1+c)M/(1-cM)$.

The case when $c$ may be taken equal to 1 in Theorem 1 deserves
special attention. Suppose $h$ is a univalent map with the
property that $h+\overline{g}$ is univalent for all $g$ with
$|g'|<|h'|$. We claim that $\Omega=h(\mathbb{D})$ must be a convex
domain. If not, it follows from [H-P] that there exists a
nonunivalent holomorphic map $\phi:\Omega\rightarrow\mathbb{C}$
with ${\rm Re}\{\phi'(w)\}>0$ for all $w\in\Omega$. It is easy to
modify $\phi$ to a mapping $\psi$ that remains nonunivalent in
$\Omega$ with ${\rm Re}\{\psi'(w)\}$ contained in a compact subset
of the right half plane $\mathbb{H}$. For example, let
$p(\zeta,t)$ be a family of univalent mappings, holomorphic in
$\zeta$ and continuous in $t$, which map $\mathbb{H}$ onto proper
subsets of $\mathbb{H}$, and which converge locally uniformly to
$p(\zeta,0)=\zeta$ as $t\rightarrow 0$. The antiderivatives of
$p(\phi'(w),t)$ appropriately chosen converge locally uniformly
$\phi$ as $t\rightarrow 0$, and  we may therefore take $\psi$ to
be the antiderivative of $p(\phi'(w),t)$ for small $t$. By scaling
the range $\psi'(\Omega)$ we may assume that it is contained in
the disk $|\zeta-1|<1$, and we can write $\psi(w)=w+\varphi(w)$,
with $|\varphi'(w)|<1$. A shown above, for suitable $\theta$, the
harmonic mapping $F(w)=w+e^{i\theta}\overline{\varphi(w)}$ will
fail to be univalent in $\Omega$, and hence
$f(z)=h(z)+\overline{g(z)}$ with $g(z)=e^{-i\theta}\varphi(h(z))$
is a nonunivalent harmonic mapping with $|\omega|<1$.

\bigskip
In our second result, the univalence of $h$ is deduced from that
of $f$.

\medskip
\noindent {\bf Theorem 2:} {\sl Let $f=h+\overline{g}$ be a
sense-preserving univalent harmonic mapping defined on
$\mathbb{D}$, and suppose that $\Omega=f(\mathbb{D})$ is linearly
connected with constant $M$. If $|\omega| < 1/(1+M)$ then $h$ is
univalent.}

\medskip
\noindent {\bf Proof:} Suppose that $h(z_1)=h(z_2)$ for distinct
points $z_1,z_2\in\mathbb{D}$. Then
$f(z_1)-f(z_1)=\overline{g(z_1)-g(z_1)}$, that is,
$$\overline{w_1-w_2}=G(w_1)-G(w_2) \, , \eqno (2) $$
where $w=f^{-1}(z)$ and $G=g\circ f^{-1}$. We will estimate the
complex derivatives $G_w=g'\cdot(f^{-1})_w$,
$G_{\overline{w}}=g'\cdot(f^{-1})_{\overline{w}}$ to show that (2)
leads to a contradiction. Differentiation of the equation
$f^{-1}(f(z))=z$ yields the relations

$$ (f^{-1})_w\cdot h'+ (f^{-1})_{\overline{w}}\cdot g' =1 \, ,$$

$$ (f^{-1})_w\cdot\overline{g'}+ (f^{-1})_{\overline{w}}\cdot\overline{h'} =0 \, ,$$
hence
$$ (f^{-1})_w = \frac{\overline{h'}}{|h'|^2-|g'|^2} \quad \,
,\quad\quad (f^{-1})_{\overline{w}} =
-\frac{\overline{g'}}{|h'|^2-|g'|^2} \, . \eqno (3) $$ It follows
that

$$ |G_w|+|G_{\overline{w}}|
=\,|g'|\frac{|h'|+|g'|}{|h'|^2-|g'|^2}\,=\,\frac{|g'|}{|h'|-|g'|}
=\,\frac{|\omega|}{1-|\omega|}\,<\, \frac{1}{M} \eqno (4) $$
because $|\omega|<1/(1+M)$ by assumption.

Let $\gamma\subset \Omega$ be a curve joining $w_1, w_2$ with
$l(\gamma)\leq M|w_1-w_2|$. Then $$|G(w_1)-G(w_2)| \leq
\int_{\gamma}\left(|G_w|+|G_{\overline{w}}|\right)|dw|
<\,\frac{l(\gamma)}{M}\, \leq \,|w_1-w_2| \, , $$ which stands in
contradiction with (2). This finishes the proof.

\bigskip

As before, one can show that when the dilatation of
$f=h+\overline{g}$ satisfies the stricter bound $|\omega| \leq c$
for some $c<1/(1+M)$, then $h(\mathbb{D})$ is a linearly connected
domain. To this effect, we regard $R=h(\mathbb{D})$ as the image
of $\Omega=f(\mathbb{D})$ under the mapping
$\phi(w)=w-\overline{G(w)}$, where $G=g\circ f^{-1}$. Let
$\zeta_i=\phi(w_i)\in R$, $i=1,2$, and let $\gamma\subset\Omega$
be a curve joining $w_1$ and $w_2$ with $l(\gamma)\leq
M|w_1-w_2|$. We will show that the curve $\Gamma = \phi(\gamma)$
satisfies
$$l(\Gamma) \leq \frac{M}{1-c(1+M)}\,|\zeta_1-\zeta_2| \, . \eqno
(5) $$ Using equations (3) we find that
$$\phi_w=1-\overline{g'\cdot (f^{-1})_{\overline{w}}}
=1+\frac{|g'|^2}{|h'|^2-|g'|^2}= \frac{|h'|^2}{|h'|^2-|g'|^2}\,
,$$
$$ \phi_{\overline{w}}=-\,\overline{g'\cdot
(f^{-1})_{w}}=-\frac{h'\overline{g'}}{|h'|^2-|g'|^2} \, ,$$
therefore
$$|\phi_w| + |\phi_{\overline{w}}| = \frac{|h'|}{|h'|-|g'|} =\,
\frac{1}{1-|\omega|} \,\leq \, \frac{1}{1-c} \, .$$ We conclude
that
$$ l(\Gamma) \leq \int_{\gamma}(|\phi_w| +
|\phi_{\overline{w}}|)|dw| \leq \frac{1}{1-c}\,l(\gamma) \leq\,
\frac{M}{1-c}\, |w_1-w_2| \, . \eqno (6) $$ On the other hand,
since $\zeta_1-\zeta_2 =w_1-w_2-\overline{G(w_1)-G(w_2)}\,$ we
have that
$$|\zeta_1-\zeta_2| \geq |w_1-w_2|-|G(w_1)-G(w_2)| \geq
|w_1-w_2|-\int_{\gamma}(|G_w|+|G_{\overline{w}}|)|dw| $$
$$ = |w_1-w_2| - \int_{\gamma}\frac{|\omega|}{1-|\omega|}|dw|
\geq \,|w_1-w_2| - \frac{c}{1-c}\,l(\gamma)$$ $$ \geq |w_1-w_2| -
\frac{Mc}{1-c}\,|w_1-w_2| = \frac{1-c(1+M)}{1-c}\,|w_1-w_2| \, .
\eqno (7)$$ The inequality (5) follows from (6) and (7).

\medskip

An interesting case of Theorem 2 is when $f=h+\overline{g}$ is a
convex harmonic mapping, that is, when $\Omega=f(\mathbb{D})$ is
linearly connected with constant 1. Then  $h$ will be univalent as
long as $|\omega| <1/2$. We do not know whether the constant $1/2$
is best possible. We also do not know, for example, if $h$ becomes
convex if $|\omega|$ is sufficiently small (and
$\Omega=f(\mathbb{D})$ is strictly convex).

In our last result, the notion of linearly connected domain is
used to describe certain classes of univalent harmonic maps $f$,
in a sense, stable.

\medskip
\noindent {\bf Theorem 3:} {\sl Let $f=h+\overline{g}$ be a
sense-preserving univalent harmonic mapping defined on
$\mathbb{D}$, and suppose that $\Omega=f(\mathbb{D})$ is linearly
connected with constant $M$. If $|\omega| <1/(1+2M)$ then
$F=h+e^{i\theta}\overline{g}$ is univalent for every $\theta$.}

\medskip
\noindent {\bf Proof:} Let $F=h+e^{i\theta}\overline{g}$ and write
$F=f+(e^{i\theta}-1)\overline{g}$. If $F$ fails to be univalent,
then for some distinct points $z_1,z_2\in\mathbb{D}$,
$$f(z_2)-f(z_1)= (1-e^{i\theta})\overline{g(z_2)-g(z_1)}\, ,$$
that is,
$$\overline{w_2-w_1}=(1-e^{-i\theta})(G(w_2)-G(w_1)) \, , \eqno (8)$$
where $w=f^{-1}(z)$ and $G=g\circ f^{-1}$. As in (4) we have that
$$ |G_w|+|G_{\overline{w}}|
 =\,\frac{|\omega|}{1-|\omega|}\,<\,\frac{1}{2M} \eqno (9) $$
because $|\omega|<1/(1+2M)$ by assumption.

Let $\gamma\subset \Omega$ be a curve joining $w_1, w_2$ with
$l(\gamma)\leq M|w_1-w_2|$. Then $$|G(w_2)-G(w_1)| \leq
\int_{\gamma}\left(|G_w|+|G_{\overline{w}}|\right)|dw|
<\,\frac{l(\gamma)}{2M}\, \leq \,\frac{|w_2-w_1|}{2} \, , $$ which
contradicts (8). This finishes the proof.

\medskip

Theorem 3 can also be deduced from Theorem 1 and the remarks to
Theorem 2, to draw a stronger conclusion. In effect, let
$f=h+\overline{g}$ satisfy the hypotheses of Theorem 3. Because
$|\omega| < c=1/(1+2M) < 1/(1+M)$, it follows from equation (5)
that $h(\mathbb{D})$ is a linearly connected domain with constant
$$ \frac{M}{1-c(1+M)} = \frac{M}{1-(1+M)/(1+2M)} =
1+2M \, . $$ Hence, we conclude from (the proof of) Theorem 1 that
any harmonic mapping of the form $h+\overline{\varphi}$ will be
univalent if $|\varphi'| < |h'|/(1+2M)$, which includes, in
particular, all rotations $\varphi=e^{i\theta}g$.

Finally, an important instance of Theorem 3 is when
$\Omega=f(\mathbb{D})$ is convex, in which case any harmonic
mapping of the form $F=h+e^{i\theta}\overline{g}$ will be
univalent if $|\omega| <1/3$. We have been unable to prove or
disprove that the constant 1/3 is sharp.

\bigskip
\noindent{REFERENCES}

\medskip
\noindent [A-B] J. Anderson and J. Becker, {\it Univalency and the
interior chord-arc condition}, preprint.

\smallskip
\noindent [C-SS] J. Clunie and T. Sheil-Small, {\it Harmonic
univalent functions}, Ann. Acad. Sci. Fenn. Ser. A.I 9 (1984),
3-25.

\smallskip
\noindent [Ch] G. Choquet, {\it Sur un type de transformation
analytique g\'en\'eralisant la repr\'esentation conforme et
d\'edinie au moyen de fonctions harmoniques}, Bull. Sci. Math. 69
(1945), 156-165.

\smallskip
\noindent [H-P] F. Herzog and G. Piranian, {\it On the univalence
of functions whose derivative has a positive real part}, Proc.
AMS, Vol. 2, No. 4 (1951), 625-633.
\smallskip

\noindent [K] H. Kneser, {\it L\"{o}sung der Aufgabe 41},
Jahresber. Deutsch. Math.-Verein 35 (1926), 123-124.
\smallskip

\noindent [P] Ch. Pommerenke, {\it Boundary Behaviour of Conformal
Maps}, Grundlehren der Math. Wiss. 299, Spinger-Verlag 1992.

\end{document}